# Euler and the German Princess

Dora E. Musielak
University of Texas at Arlington
dmusielak@uta.edu

**Résumé**. In 1760, mathematician Leonhard Euler began to write beautiful *Letters to a German Princess on Diverse Subjects of Physics and Philosophy*. Much has been written about Euler and his work, but we wonder, *who* was the princess? How did she become involved with the greatest mathematician of her time? The princess was a fifteen-year old named Friederike Charlotte von Brandenburg-Schwedt. In this article we explore her story.

**Résumé**. En el año 1760, el matemático Leonardo Euler empezó a escribir las bellísimas *Letters to a German Princess on Diverse Subjects of Physics and Philosophy*. Mucho se ha escrito acerca de Euler y su obra, pero nosotros nos preguntamos, ¿Quien era esa princesa? ¿Cómo se asoció con el genio matemático más celebrado de su tiempo? La princesa de quince años era Friederike Charlotte von Brandenburg-Schwedt. En este articulo exploramos su historia.

## 1.     The Princess

Princess Friederike Charlotte belonged to the Brandenburg-Schwedt line of the Prussian royal family. She was born on 18 August of 1745 in Schwedt, formerly a Prussia town in northeastern Brandenburg, now Germany. With the status of a Große kreisangehörige Stadt, it is now the largest town of the Uckermark district, located near the Oder River, about 100 km north of Berlin and close to the border with Poland. Friederike Charlotte was the firstborn daughter of Friederich Heinrich von Brandenburg-Schwedt and his wife Leopoldine Marie von Anhalt-Dessau, a daughter of Prince Leopold I of Anhalt-Dessau. Friederike Charlotte had a younger sister named Louise Henriette born in 1750.

The princess's father, Frederick Henry, was the last Margrave of Brandenburg-Schwedt, owner of the Prussian secundogeniture,[1] ruling from 1771 to 1788. Margrave was the medieval title for the military commander assigned to maintain the defense of one of the border provinces of the Holy Roman Empire. Married since February 1739, Frederick Henry and Leopoldine Marie were very unhappy. After the birth of their two daughters, the couple quarreled often and violently and soon they separated. Friederike's mother went to live in Kołobrzeg, a palace near the Parsęta River on the south coast of the Baltic Sea. Today Kołobrzeg is a city in the West Pomeranian Voivodeship in northwestern Poland with some 50,000 inhabitants.

---

[1] A secundogeniture was a dependent territory given to a younger son of a princely house and his descendants, creating a cadet branch.

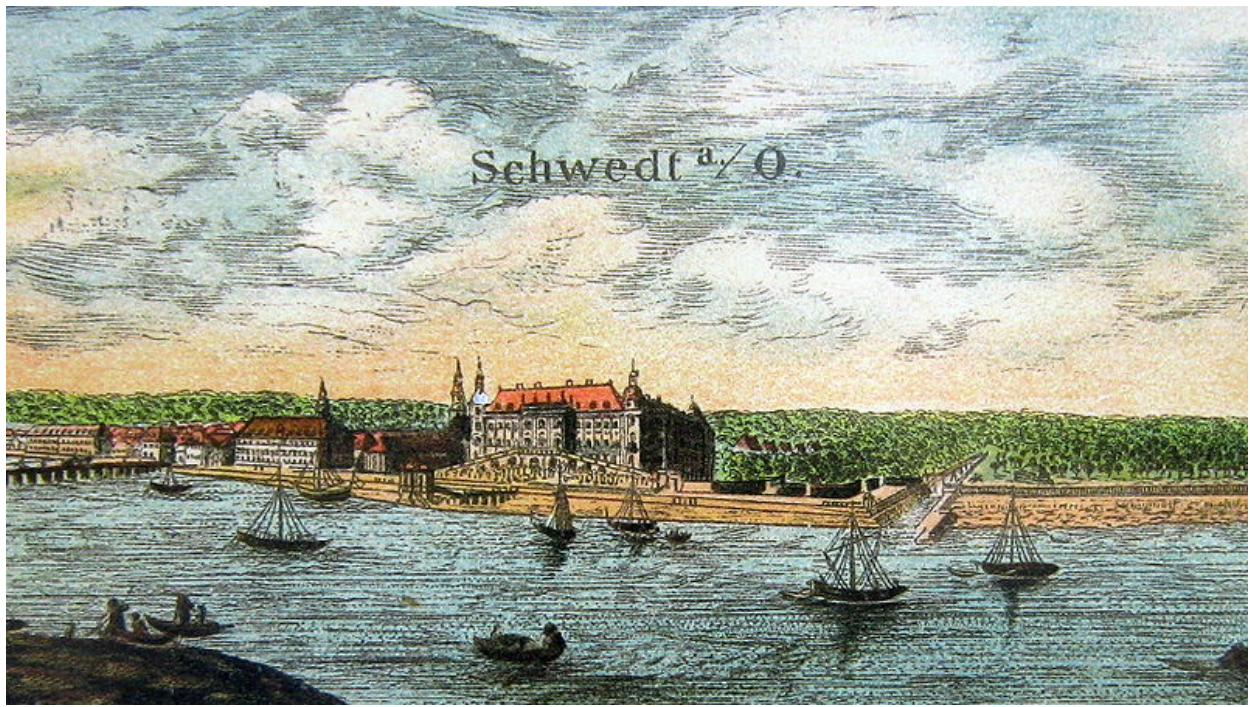

Schwedt Castle overlooking the Oder River, 1669. It was destroyed in 1945.

As it was customary in those days, Friederike Charlotte and her sister were entrusted to the care of nannies and were perhaps taught by governesses at one or both of their parents' castles. It is likely that Friederike Charlotte received her basic education at the Herford Abbey—girls of the nobility were educated at convents to prepare them for marriage— and only visited her parents on special occasions. Her grandmother, Johanna Charlotte von Anhalt-Dessau was Princess-Abbess at Herford Abbey. Johanna became the new abbess on 10 October 1729, but not until 1735 she took up permanent residence in Herford and appointed Hedwig Sophie of Schleswig-Holstein-Gottorp [2] as her coadjutor. By age ten, Friederike Charlotte was named coadjutor (*Koadjutorin*) to Abbess Hedwig Sophie who had succeeded her grandmother in 1750. Coadjutor meant the young princess was designated as the successor of Abbess Hedwig Sophie.

The residents at the Abbey were the daughters of noble aristocratic families. Many had arrived as little girls to be educated by the nuns. Others were young women with no prospects for marriage or princesses who were groomed to rule as princess abbess. For example, at nineteen Princess Therese Natalie von Brunswick-Wolfenbüttel-Bevern became a collegiate lady in Herford Abbey because efforts to marry her to a French prince failed. Shortly after, it was decided that she would succeed Elisabeth of Saxe-Meiningen as Abbess of Gandersheim Abbey.

---

[2] Hedwig Sophie Auguste von Schleswig-Holstein was born on 9 October 1705 at Gottorp Palace, Schleswig, Schleswig-Flensburg, Schleswig-Holstein, Germany. She was Prioress of Quedlinburg in 1728, and Canoness of Herford. in 1745. She became Abbess of Herford. in 1750. She died on 13 October 1764 at Herford, Detmold, North Rhine-Westphalia, Germany, at age 59. From Detlev Schwennicke, *Europaische Stammtafeln*, New Series, Vol. I/3, Tafel 293.

How Euler become acquainted with Princess Charlotte is not immediately clear. Let us explore some historical facts that may provide some clues to their relationship.

## 2.     The Mathematician Leonhard Euler

When Frederick II ascended the Prussian throne in 1740, Leonhard Euler was thirty-three years old, already known throughout Europe as a *célèbre mathématicien*. In addition to having published numerous articles and books in several branches of mathematics, Euler, who was academician professor at the Saint Petersburg Academy of Sciences in Russia, had just won the prestigious prize of the French Academy of Sciences. He was a coveted scholar. At the same time, one of the first initiatives of the Prussian King was to revitalize the Academy of Sciences in Berlin, which had languished for the past decades. To give it the prestige of the French Academy, Frederick II sought out the most prominent men of science, extending a special invitation to Euler to join his Royal Prussian Academy.

Euler accepted, arriving with his family in Berlin in the summer of 1741. By the time princess Friederike Charlotte was a teen, Euler had amassed a large scientific legacy and was undisputable the greatest mathematician of his time. He was busy. In addition to his enormous scholarly work, writing numerous articles and several major treatises, Euler was entrusted with many duties at the Academy—he had to oversee the astronomical observatory and botanical gardens, deal with personnel matters, attend to financial affairs, notably the sale of almanacs, which constituted the major source of income for the Berlin Academy. The king also charged Euler with engineering projects, such as correcting the level of the Finow Canal in 1749. At that time, he also supervised the work on pumps and pipes of the hydraulic system at Sans Souci, the royal summer residence.[3]

Euler was known not just among the scientific circles but also among the members of the Prussian court. It is possible that Friederick Henry von Brandenburg-Schwedt, the princess's father, became acquainted with Euler in 1755 when he purchased the *Prinzessinnenpalais*, a stately mansion in the center of Berlin where the Prussian nobility resided. He and Euler shared an interest in music and could have met at some musical concert. Music was prominent in Prussia. Being a musician himself, Frederick II's court was a center of German musical life and included prominent musicians such as C. P. E. Bach, second son of the famous J. S. Bach.

Friederick Henry was a patron of the arts. As soon as he took the title of Margrave of Brandenburg-Schwedt in 1771, he developed Schwedt into a cultural center. He had an operetta theatre for four hundred people built in the orangery of his castle. Thanks to him, Schwedt became one of the first theatre towns in Germany.

Even though we don't know how or when Euler met the young princess in person, he wrote her extensive letters between April 1760 and 1762, addressing a wide range of topics dealing with science and philosophy. The scope of the subjects he covered suggests that Euler's intention was to tutor her in preparation for something beyond the casual curiosity of a young girl. Who asked Euler to tutor the princess and why?

---

[3] A P Youschkevitch, Biography in Dictionary of Scientific Biography (New York 1970-1990).

It is speculated that Euler visited the future Margrave in his Berlin castle with the objective of tutoring his daughters.[4] In 1782, Condorcet implied that the princess wished to receive from Euler a few lessons in physics.[5] However, by 1760 Charlotte was already coadjutor at Herford Abbey, and it is much more likely that the lessons in the letters were intended to prepare her to rule as Princess Abbess.

The Imperial Immediacies (*Reichsfreiheit* or *Reichsunmittelbarkeit*) held a privileged feudal and political status and had the right to collect taxes and tolls and held juridical rights themselves. They functioned as landowners, regional administrators, and exercised the right of the lower courts. Thus, the teenager princess had to be fully prepared to exercise those functions. Euler, the supreme scholar at the Prussian court, was the most indicated teacher to guide her, helping to develop her intellect. He could also advise the princess in philosophical and theological matters. Euler, who was a Calvinist, was much admired for his moral and religious character. He never abandoned the religious duties to which he had been educated. It was known that, as long as he preserved his eyesight, he assembled his entire family every evening and read a chapter of the Bible, which he accompanied with an exhortation.

When Euler began to write to Princess Friederike Charlotte, he was fifty-three years old and partially blind—he had lost vision in his right eye since he was thirty-one during a period of intense work at the St. Petersburg Academy. In 1760, Euler was at the peak of his career as Director of the Mathematics Section of the Prussian Academy of Sciences in Berlin. It is rather amazing that Euler would have had time to instruct anyone. Yet, he did.

On the 19th of April 1760, Euler penned the first letter to the fifteen-year-old princess. He addressed her as *Votre Altesse*, Your Highness. Over the course of two years he wrote two hundred-thirty four letters. Euler composed them on Saturdays and Tuesdays, sometimes writing two or three different letters on the same day to lecture on different topics, laying out the basics of a dozen disciplines. He taught Princess Friederike Charlotte almost everything —physical science, astronomy, music, logic, theology, and philosophy. That was the Age of Enlightenment, and Euler wrote in French, which was the language of the Prussian court. He explained many important topics of natural science and provided their philosophical background in clear and concise manner that would be understandable to the young girl.

The year 1760 was a troublesome time for Prussia due to the increased foreign hostilities of the Seven Years War (1756-1763). The queen, members of the royal family, the archives, and the principal ministries had left for Magdeburg (over 150 km SW from Berlin) after the disaster at Kunersdorf in 1759. The Raid on Berlin took place in October 1760 when Austrian and Russian forces occupied the Prussian capital for several days, left vulnerable by Frederick's decision to concentrate his forces in Silesia. The occupiers ransacked many areas of Berlin, and several royal palaces were burnt. They demanded ransom money from the city's governor, but then a rumor that King Frederick was marching to the rescue of Berlin with his superior forces prompted the occupying commanders to withdraw. Having completed their major objectives, they left Berlin on the 12th of October.

---

[4] Calinger, R., *Euler's "Letters to a Princess of Germany" As an Expression of his Mature Scientific Outlook*. Archive for History of Exact Sciences, Vol. 15, No. 3 (23.VII.1976), pp. 211-233. Published by Springer.

[5] Condorcet wrote: Madame la princesse d'Anhalt-Dessau, nièce du roi de Prusse, voulut recevoir de lui quelques leçons de physique … p. 17. Moreover, Condorcet misidentified the princess.

Euler's country estate in Charlottenburg was mistakenly looted during the occupation. Euler wrote to his friend Müller: "There, almost everything has been destroyed and devastated, because I have received the security guard that His Lordship, Count Chernishev has graciously granted me, too late, for as soon as our army had retreated, Charlottenburg was immediately surrendered to the Cossacks. I thereby have lost 4 horses, 12 cows, a good many small animals, a lot of oats and hay, and in the house all the furniture is ruined."[6] Euler was later indemnified.

Now, if he began tutoring the princess at her Berlin residence, then why did Euler switch to writing the lessons? It is most likely that by the spring of 1760, Friedeike Charlotte was no longer in Berlin and Euler felt compelled to continue the instruction by correspondence. In fact, the opening paragraph of the first letter suggests that Euler had lost hope of giving her instruction in person (see below). In letter 47 he wrote "Were you to draw a straight line from any point, in your apartment at Berlin, to a given point in your apartment at Magdeburg,..."[7] suggesting that Friederike was in Magdeburg. It is highly likely since the queen and her court had moved there the previous year.

Four years after concluding the instruction of the princess, Euler accepted the invitation of the Russian empress Catherine II to return to the Academy of Sciences in St. Petersburg. Soon after, Euler became almost entirely blind after an illness. Someone in Russia must have read the copies of the letters Euler had written to the princess and discovered the scientific and philosophical richness they contained. The scope and depth of the topics Euler treated made the collection of letters as a unique encyclopedia. The empress encouraged the publication of the letters, which would serve to make science accessible to a wide range of readers. She was right. Published in 1768 in the original French language that Euler wrote, the *Lettres à une Princesse d'Allemagne sur divers sujets de Physique & de Philosophie*[8] became an instant success. By the end of the eighteenth century, the three-volume book had been translated into almost every European language, and had gone through several dozen printings.[9]

The first edition of the *Lettres* published in Paris in 1787 included *Eloge de M. Euler*, a thirty-six page obituary written by the marquis de Condorcet that provided the reader with a biographical sketch and highlights of Euler's career.[10] Condorcet was Secretary of the French Académie des Sciences since 1777. Even though Euler wrote the letters in French, it is evident that Condorcet made some editorial changes because the text differs from the original. Perhaps Condorcet wished to impart a more elegant Parisian style to Euler's original manuscript. This is evident from the first letter, which Euler wrote on 19 April 1760:

"Madame, Comme l'espérance de pouvoir continuer à votre Altesse mes instructions dans la Géométrie semble de nouveau être reculée, ce qui me cause un très sensible chagrin, je souhaiterais y pouvoir suppléer par écrit, autant que la nature des objets le permet." (1770)

Condorcet edited this paragraph as follows:

"Madame, Comme l'espérance de continuer à votre Altesse mes instructions dans la Géométrie semble s'éloigner de plus en plus, ce qui me cause un chagrin très sensible, je souhaiterais pouvoir y suppléer par écrit, autant que la nature des objets peut le permettre." (1782)

In free English translation we read:

"Madame, as the hope to continue communicating in person to Your Highness my lessons in Geometry seem more remote, for which I am sorry, I feel compelled to provide you personal instruction by writing, as far as the nature of the subjects allow it." (1770)

"Madame, as hope to continue giving Your Highness my lessons in geometry seems to recede away more and more, which causes me much grief, I wish to compensate for them by writing, as much as the nature of the topics can allow it." (1782)

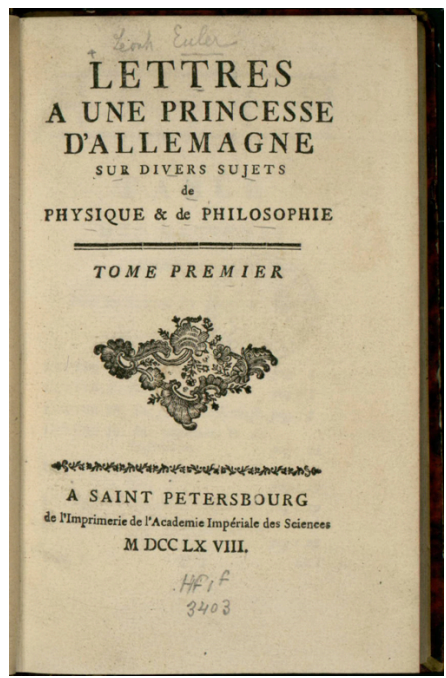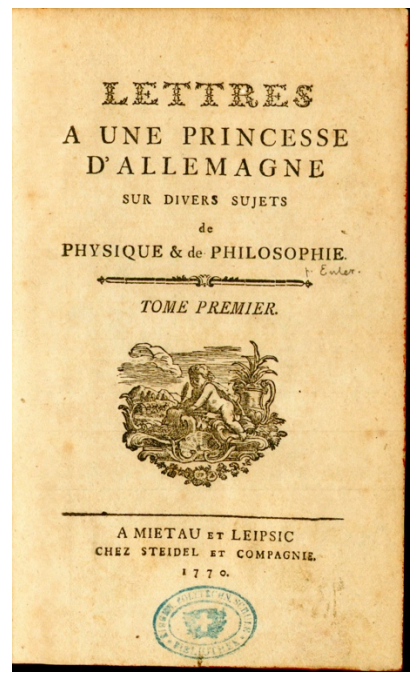

Frontispieces of *Lettres à une Princesse d'Allemagne sur divers sujets de Physique & de Philosophie* (1768, 1770).

### 3. Excerpts from Selected Letters[11]

The nature and the scope of the topics that Euler addressed in those letters to Princess Friederick Charlotte are of interest even today. He wrote about gravity, astronomy and the laws of motion, the nature of sound and light, electricity and magnetism, always using a clear manner and using an effective pedagogical style without invoking equations or formulas. Euler thoroughly explained all the topics that should be included in an introductory textbook of physics of the time, including beautiful illustrations when needed.

Let us read a few excerpts from the *Lettres*, just enough to understand what the young princess learned from Euler. From the first letter one realizes that the reader will be exposed to an introductory course to physical science. Euler intended for his young correspondent to become familiar with units so she could better understand the topics that would follow. He explained the concept of magnitude (*la grandeur*), giving examples from "the smallest to the greatest extensions of matter discoverable in the world." Euler gave special emphasis to units of length such as feet and miles, with examples to compare them, saying that is more clear to quote the distance from Berlin to Magdeburg as 18 miles rather than 432,000 feet.

In addition to explaining the theories of light and gravity expounded by Newton and Descartes, Euler discussed optics, electricity and magnetism. At least seventeen letters deal with electricity (Letters 138-154) and nineteen letters cover magnetism (Letters 169-187), suggesting that these topics were popular at the time. In fact, Euler had won in 1744 and 1743 two prizes from the French Academy of Sciences for his work on magnetism. Euler would be pleased to know that exactly a hundred years later, James Clerk Maxwell would synthetize the extensive knowledge available at the time in electromagnetism and light into four concise, and rather beautiful mathematical equations.

In other letters, Euler also spent a great deal explaining issues related to the functioning of scientific devices such as the telescope, thermometers, and microscopes. And perhaps because of his own vision problems, Euler taught her about vision and glasses, including detailed descriptions the structure of the human eye, the eye of an animal, and the artificial eye, or *camera obscura* (Letters 40-44).

Euler wrote several letters to teach the princess about the composition of air, its gravity, rarefaction and expansion, and changes in the atmosphere caused by heat, cold, and altitude. In letter 15, he explained how air changes when exposed to heat and cold: "Heat and cold produce the same effect on air, as on every other body. Air is rarified by heat, and condensed by cold. However, from what I have explained to Your Highness, a certain amount of this air is not determined to occupy only a certain space, as all other bodies are; but by its nature, air always tends to dilate, and actually does expand itself, as long as it meets no obstacle. This property of air is denominated elasticity" Euler went on, concluding that the elasticity of air is in exact proportion to its density.

Euler's words are reminiscent of Charles's law, which describes how gases tend to expand when heated. Let's keep in mind that the ideal gas law (the equation of state of

---

a hypothetical ideal gas that we use as a good approximation to the behaviour of many gases) was not known at the time. It was not until 1834 when French engineer Émile Clapeyron combined Charles law and Boyle's law that the ideal gas law was established in the form $P = \rho RT$, establishing the relationship between pressure, density, and temperature, a formula which is valid for gases over a wide range of conditions.

Letter 17 addressed light and the theories of Descartes and Newton. Euler began with the rays of the Sun, "which are the focus of all the heat and light that we enjoy." Then he asked, what is the nature of these rays? He thought that it is one of the most important inquiries in physics, because "from it an infinite number of phenomena is derived." Euler was correct in pondering the nature of the solar rays. It's rather amazing to think that for most of human history, we had very little understanding about our own star. We didn't know what the Sun was made of, how it formed, or how it produced energy. In Euler's time astronomers didn't know how big the Sun was, and its distance from the Earth was rather approximated. Regarding this distance, Euler was far off. In the previous letter he wrote that it was about "trente Millions de Milles" (thirty million of miles). He used German miles,[12] so his value corresponds to approximately 272.4 million kilometers. In reality, we orbit the Sun at an average distance of about 150 million kilometers. Moving in an elliptical path, at its closest point the Earth gets to 147 million km from the Sun, and at its most distant point, it is 152 million km.

Of course Euler knew that the distance between the Earth and the Sun was estimated, starting in antiquity. In 1761, the transit of Venus would give astronomers the opportunity to measure accurately the size of the Solar System. Thanks to a rare celestial alignment, Venus would pass in front of the Sun, taking about six hours to cross it. By recording the times of the start and end of the event from widely separated locations around the Earth, astronomers could use trigonometry to calculate the distance to Venus and then to the Sun. Having that, they could use Kepler's laws of planetary motion to calculate the orbits of all the planets out to Saturn, the farthest known planet. Due to many unforeseen challenges, it was not until 1771 that astronomers were able to make the final calculation of the distance to the Sun. The combined results from all the various missions were within about 4% of the modern accepted value. At the next pair of transits, in 1874 and 1882, the accuracy was improved to 1%. In the twentieth century, radar experiments with radio telescopes improved those results and led to the definition of the Astronomical Unit.

Let's go back to 1760. Regarding the nature of light, Euler inquired: "Do certain portions, inconceivably small, of the Sun itself, or of its turbulence, come down to us?" He wondered whether the transmission of the solar rays was similar to the sound of a bell. Then Euler wrote about Descartes, who believed that the whole universe was filled with matter composed of small globules, which Descartes called the second element. That theory supposed that the Sun was in a state of agitation, transmitting it to the globules, and that "they again communicate their motion in an instant to every part of the universe." However, Euler noted, the solar rays (light) do not reach us instantaneously bur rather they take eight minutes to travel through the Sun-Earth distance, and thus the theory of Descartes had been discarded.

Indeed, Danish astronomer Ole Christensen Rømer made the first quantitative measurements of the speed of light in 1676. He used the inequality of Jupiter's moons.

---

[12] One German mile was about 9.08 km.

Rømer observed that times between eclipses (particularly those of Io) got shorter as Earth approached Jupiter, and longer as Earth moved farther away. Rømer's conclusion was that the speed of light is finite. This was not fully accepted until 1727 when James Bradley made measurements of the so-called aberration of light.

Euler wrote about Newton's theory, which maintained that the luminous rays are separated from the body of the Sun, and "the particles of light thence emitted with that inconceivable velocity which brings them down to us in eight minutes.... This opinion is called the system of emanation: it being imagined that rays emanate from the Sun and other luminous bodies, as water emanates or springs from a fountain." To Euler, this theory was bold and irreconcilable to reason. He explained why he thought Newton's system of emanation theory was unreasonable and he did not accept it, even if many philosophers did. Euler ended the letter by saying that he was "too little a philosopher to adopt the opinion in question." (Pour moi je suis trop peu Philosophe pour embrasser ce sentiment.)[13]

Clearly, Euler refuted Newton's idea, and in the next letter he expanded on the reasons why. The question whether light is comprised of streams of particles (corpuscles) or made up of waves is a very old one. Christian Huygens believed that light behaved as waves that spread out from the source that generates the light, and thus each color is a different wavelength. On the other hand, Newton developed a corpuscular theory that treats light as being composed of tiny particles. We use Newton's theory to describe reflection, and while it can explain some physical phenomena, it cannot explain others.

During Euler's time, scholars and philosophers did not consider that light could behave as both waves and corpuscles. Early in the nineteenth century experiments by Thomas Young suggested that light is a wave motion. In 1803, the results of his studies on diffraction and interference of light supported Huygens and thus opposed Newton's theory. Years later, the experimental research of French engineer Augustin-Jean Fresnel showed that light is a transverse wave; he used the so-called Fresnel's integrals to validate his theory. Then in 1865, James Clerk Maxwell demonstrated that electric and magnetic fields travel through space as waves moving at the speed of light. Maxwell proposed that light is in fact undulations in the same medium that is the cause of electric and magnetic phenomena. Maxwell's theory is an improvement over the wave theory in that it explains how light is generated. Visible light is only one type of electromagnetic wave. This theory is very mathematical intensive. Finally, at the turn of the twentieth century, scientists developed the Quantum Theory of Light which describes light as (in some sense) both a particle and a wave, and (in another sense), as a phenomenon which is neither a particle nor a wave. This theory is more mathematical than Maxwell's equations. It is similar to the corpuscular theory of Newton, except that instead of describing light as particles, light is composed of packets of energy or photons. The energy of the photon determines the color. Thus, Euler was ahead of his time sharing with a young student his uncertainty about the theory of light, something that took centuries to be fully established.

In letter 19, Euler wrote "You have seen that the system of the emanation of the rays of light...." and near the end of the letter he adds, "Using certain artifices, we can cause mercury to emit light, as Your Highness remembers seeing."[14] This suggests that

---

[13] Euler, L., Letter XVII, p. 69
[14] Dominic Klyve, Personal Correspondence, June 2013.

Euler and the princess witnessed the same experimental demonstration. Euler explained the nature of sound and then introduced its velocity of propagation in order to address consonance and dissonance and the "pleasure derived from fine music." Euler was clearly an authority in the subject of sound. His first published memoir was *Dissertatio physica de sono,* written when he was twenty-years old.

Now imagine someone explains Newton's *Principia* to you in a very lucid and concise manner. If there were a teacher who could explain the laws of motion and universal gravitation, without using a single equation or formula, and who would teach you how the planets move around the Sun and why the Moon doesn't fall, that teacher would be Euler. In the course of twenty-four letters, Euler taught the teenage princess about gravity, the laws of motion, the reason why the Moon keeps going around our planet, and how it exerts a force that causes the tides. He even discussed the different philosophical ideas about gravitation, and he also taught her about the Solar System with its planets and moons and how the combined forces of attraction make their motion so difficult to model.

Euler introduced the concept of *gravity*, starting Letter 45 by writing that all bodies, solid and fluid, fall downward, when they are not supported. "I hold a stone in my hand; if I let it go, it falls to the ground, and would fall still farther, were there an aperture in the Earth. The paper on which I write would fall to the floor, were it not supported by the table. The same law applies to every body. There is not one that would not fall to the ground, if it were not supported, or stopped. The cause of this phenomenon, or of this propensity of all bodies, is named their gravity or their weight." He added that when bodies are said to be heavy, or possess gravity, we mean, that they have a propensity to fall downward, and actually would fall, if we remove what before supported them. In the next five letters, Euler wrote about gravity extensively, including physical and philosophical discussions to address the true direction and action of gravity relatively to the Earth. Then on Letter 51 Euler was ready to present in detail the action of gravity on the Moon, explaining why it doesn't fall.

Euler drew a similarity with a body in motion, using the example of a cannonball shut from the top of a mountain. These are Euler's words, reproduced in the original orthography, and then in free English translation:

Un boulet de canon tiré selon une direction horizontale, ne parvient à la terre que sort loin ; & si l'on le tiroit sur une haute montagne il parcourroit peut-être plusieurs milles avant que d'arriver à la terre. Qu'on hausse encore davantage le canon, & qu'on augmente la force de la poudre, & le boulet alors sera porté beaucoup plus loin. On pourroit pousser la chose si loin, que le boulet ne tomberoit que chez nos Antipodes,[15] & en la poussant encore plus loin il pourroit arriver que le boulet ne tomberoit plus du tout, mais qu'il retourneroit à l'endroit où il a été tiré, & seroit ainfi un nouveau tour du monde ; ce seroit une petite lune qui seroit ses révolutions de même que la véritable autour de la terre. Que V. A. daigne à présent réflechir sur la grande hauteur où la lune se trouve, & sur la prodigieuse vîtesse dont elle est portée. Elle ne sera plus surprise alors que la lune ne tombe pas à terre, quoiqu'elle soit poussée par la gravité vers son centre. Une autre réflexion mettra cela dans un plus grand jour. Nous n'avons qu'à bien considéter le chemin qu'une pierre jettée obliquement, ou un boulet de canon, décrit. Le chemin est toujours une ligne courbe, telle que représentes la figure ci-jointe.

---

[15] Antipodes are places diametrically opposite each other on the globe

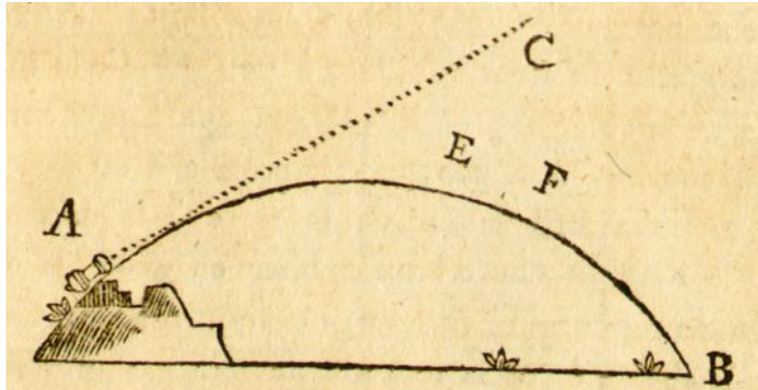

Euler's illustration to explain the force of gravity on a cannonball, Letter LI, p. 219.

"A cannon ball fired in a horizontal direction does not come down to the Earth till it has got to a considerable distance; and if it were fired from the top of a high mountain, it might, perhaps, fly several miles before it reached the ground. If the the cannon is further elevated, and we increase the strength of the gunpowder, the ball will be carried much farther. One may propel it so far, that the bullet would reach the antipodes; and propelling it still further the bullet would not fall at all, but return to the place where it was shot off, and would thus perform a new circuit around the Earth. It would thus be a little moon, making it's revolutions round the Earth like the real moon. As Your Highness please reflect on the height of the Moon, and on the prodigious velocity with which it moves It won't be more surprising that the Moon does not fall to the ground, even though it is pushed by gravity toward it's center. There is another idea which will place this in a clearer light. We have only to consider the path described by a stone thrown obliquely, or a cannonball just described. The path is always a curved line, such as represented in the annexed figure."

After that explanation, Euler emphasized that the Moon does not move in a straight line, and that its path must necessarily be curved, such as a circle around the Earth with a radius equal to the distance to the Moon. Euler knew a lot about the motion of the Moon since he had developed the first lunar theory a few years earlier. Now we know that the Moon's orbit around the Earth is elliptical, with a substantial eccentricity $e$ = 0.0549 compared with Earth's, which is 0.017. The lunar orbit is tilted about 5° to the plane of the Earth's orbit. On average, the Moon is at a distance of about 385000 km from the centre of the Earth, which corresponds to about 60 Earth radii. The Moon moves with a mean orbital velocity of 1.023 km/s.

After having taught the princess the law of universal gravitation, he was ready to teach her astronomy, the branch of science that he helped become a rigorous discipline.[16] With a combination of theoretical and practical advances, astronomers in the mid 1700 began to solve most of the outstanding problems that had been raised by Newtonian mechanics. New instruments allowed the cataloging of stars and new theories were attempting to explain the anomalies observed in the motion of celestial bodies across the sky. Astronomers discovered new phenomena such as aberration, nutation, nebulae, and astronomical expeditions were determining the size and shape of the Earth. At the same time, Euler's theoretical work aimed to answer the outstanding questions raised by Newton's theory of gravity. He had won four prizes (in 1740, 1748, 1752, and 1756) for his work in astronomy. No question about it, Euler was at the head of the group of mathematicians developing new theories in an effort to understand and explain the motion of the celestial bodies. Despite his partial blindness, Euler could visualize how the planets orbit the Sun, and that vision combined with his mathematical genius and mastery of the laws of mechanics made Euler the perfect teacher of astronomy.

In Lettre 59, Euler described our Solar System (*système du monde*), explaining the movement of the planets (*corps célestes*) around the Sun and the forces that cause their motion. He used a very beautiful illustration of the Solar System to complement his explanations; the schematic shows the six known planets: Mercury, Venus, Earth, Mars, Jupiter, Saturn, and also the known moons—our Moon plus four of Jupiter and five of Saturn. He represented the highly elliptical orbit of comets. In Lettre 60 Euler clarified that "the orbits of the planets are not circles, as the figure seemed to indicate, but rather somewhat elliptical (*ovale*), some more, some less, so that these planets are sometimes nearer to the Sun, sometimes farther off." He noted that the orbit of Venus is almost a perfect circle; but those of the other planets are more or less elongated. It is true; the eccentricity of Venus's orbit is 0.0068, while that of the Earth is 0.0167. Mercury has the most elliptic orbit with an eccentricity of 0.2056.

Euler asked the princess to observe that "the fixed stars are bodies fully similar to the Sun, and glow in the same way," and that they are very far, at "prodigious distances," adding that he already had the honour to say to her majesty that "the fixed star which is closest to us, is at least 400000 farther from us than the Sun." Moreover, he wrote, "Every fixed star seems to be intended to warm and to illuminate a number of opaque bodies, similar to our Earth, and also inhabited without doubt, which lie in its vicinity, but that we see as points of light because of their prodigious distance." Euler was predicting the existence of extrasolar planets!

In Euler's words, "... astronomy is the science that depends on an exact knowledge of the movement of all celestial bodies, in order to be able to determine, for every instant of time, both past and future, where each celestial body should be, and in what place in the sky it must appear, whether observed from the Earth or from another place whatever of the universe. The science that deals with motion in general is named *mechanics* or dynamics. Its purpose is to determine the movement of all bodies, when they are driven by any forces. This science is one of the principal branches of mathematics, and those who apply to it exert all their efforts to bring the mechanics to its highest degree of perfection."[17]

---

[16] Musielak, D., *Euler: Genius Blind Astronomer Mathematician*. Freiburg, Germany. June 2014. In this monograph I review Euler's most important contributions to astronomy.

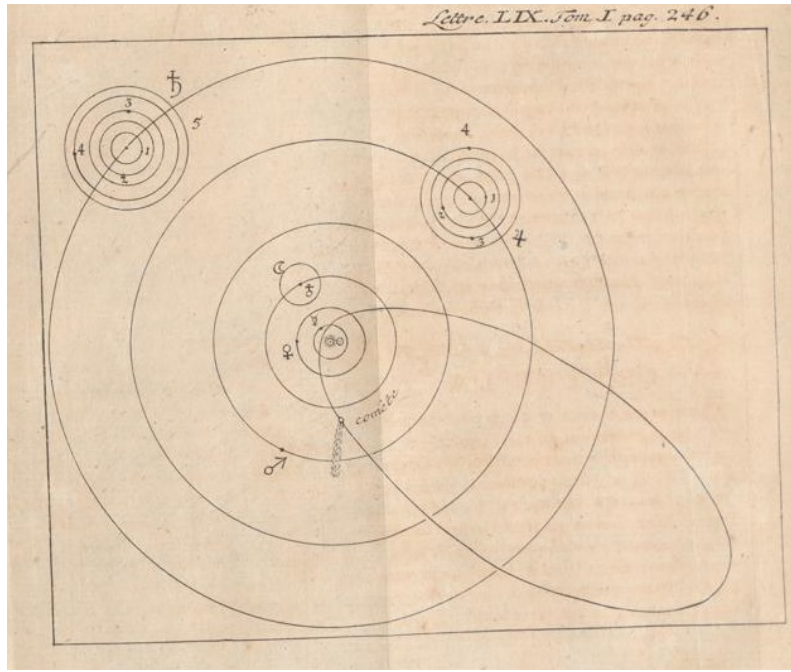

The Solar System as it was known in 1760. From Letter LIX, p. 247  (1770)

Euler admitted that the subject is rather complex. He believed that "the greatest challenge is in the number of forces which act on each celestial body, being pushed or drawn towards others. If each body were attracted to just one other body, there would be no difficulty; the great English mathematician, Mr Newton, who died in 1728, was the first to determine the motion of two bodies that attract each other, according to the law which I had the honour of speaking to Your Highness. Due to this law, if the Earth was drawn only to the Sun, we would know perfectly how to determine the movement of the Earth, without need of more research. The same thing would apply to the other planets, Saturn, Jupiter, Mars, Venus, and Mercury, if these bodies were attracted only by the Sun. Since the Earth is attracted not only by the Sun, but also by all the other heavenly bodies, the problem becomes infinitely more complicated and difficult, because the diversity of the forces that must be considered. You may neglect, however, the forces with which the Earth is attracted towards the fixed stars, because although their masses are enormous they are at very far distances that the forces they exert on the Earth is considered very small." This was a very lucid explanation of the universal law of gravity.

Seemingly oblivious to the foreign occupation during the Raid on Berlin on 9 and 12 of October, Euler continued writing; this time the topic was on the ocean tides caused by the force of attraction of the Sun and the Moon. Euler knew quite a lot about tides, having won a prize from the French Academy of Sciences in 1739 for an explanation of

---

the ebb and flow of the sea (*le flux et le reflux de la mer*). For this topic, he wrote to the princess four separate letters.

Having already treated the attractive gravitational forces to which the Earth is subjected by the Sun and the Moon, Euler proceeded to describe the ocean tides, or as he put it, the *flux and reflux of the sea.* He wrote: "We have just seen that the different parts of the sea are attracted unequally by the Moon, according as they are unequally distant from her centre; the sea must, therefore, be agitated by the force of the Moon, which, continually changing her situation, with respect to the Earth, and performing a revolution round it in about twenty-four hours and three quarters, makes the sea undergo the same changes, and presents the same phenomena in the same period of twenty-four hours and three quarters; the flux and reflux must, therefore, be retarded from one day to another three quarters of an hour, which is confirmed by constant experience. It now remains that we know, how the alternate elevation and depression of the sea, which succeed each other after an interval of fix hours and eleven minutes, result from the inequality of the powers of the Moon."

Euler also discussed that having discarded the theory of Descartes to explain the tides by the pressure of the Moon on our atmosphere, it was more reasonable to look for the cause in the attraction that the Moon exerts on the Earth, and consequently also upon the sea. He emphasized that the attractive force of all celestial bodies have been sufficiently established by so many other phenomena, as he had already explained to the princess, there was no doubt that the flux and reflux of the sea must be an effect of that force. Euler added: "As soon as we establish that the Moon, as all other celestial bodies, has the strength to attract all bodies in the direct ratio of their mass, and in the inverse ratio of the square of their distance, it is easy to see that its action extends to the sea, a body of fluid, as Your Highness must frequently have observed, that the slightest force is able to shake a fluid. All that remains is to inquire whether the attractive force of the Moon, as we assume it, is actually capable of producing in the sea the agitation that we know by the name of ebb and flow."[18]

I do not know if Princess Friederike Charlotte had ever seen the ocean, but if she hadn't Euler's description was so vivid that one could understand its rhythmic up and down motion. Euler himself had sailed through the Baltic Sea, and he noted that the tidal motion "in seas bounded and confined, such as the Baltic, and the Mediterranean, it is much less considerable ..."

Euler taught the princess the concepts of uniform, accelerated, and retarded motion. In letter LXXIII, Euler wrote Newton's laws of motion as follows: This principle of motion is commonly expressed in the two following propositions: First; *A body once at rest will remain eternally at rest, unless it be put in motion by some external or foreign cause*: Secondly; *A body once in motion will preserve it eternally, in the fame direction, and with the same velocity; or will proceed with an uniform motion, in a straight line, unless it is disturbed by some external, or foreign cause.* In these two proportions consists the foundation of the whole science of motion, called *mechanics*.

Euler was in fact the utmost expert in mechanics. At twenty-nine, he had published the first volume of *Mechanica*,[19] a book in which he restated Newton's laws

---

and discovered new features of motion. Euler's approach was superior to Newton's, as he applied the full analytic power of calculus to dynamics. In this book, Euler recast Newton's theories of motion in a more sophisticated mathematical language. In fact, Euler was the first scholar to write Newton's second law in a differential form. The familiar formula of the second law ($F = ma$) first appeared in a memoir that Euler published in 1750. With his work, Euler laid the foundations of analytical mechanics and outlined a program of studies embracing every branch of science, involving a systematic application of analysis. It is a pity that he did not teach calculus to Friederike Charlotte. But then again, she was destined for a different career, as we shall see.

Apart from physical science, Euler discussed philosophical and theological issues. Theology was one of his favorite subjects. In the last few letters of the first volume, Euler considered the action of the soul upon the body, and of the body upon the soul, and he talked about moral and physical evil. He gave his opinion on the real destination of Man, believing in the usefulness and necessity of adversity, and he added his thoughts on true happiness. Euler's letter dated 31 March 1761 deals with "the true foundation of human knowledge."[20]

## 4.    Princess Abbess

When Abbess Hedwig died, on 13 October 1764, Friederike Charlotte became the Princess Abbess of Herford, taking the title *Friederica Charlotta Leopoldina Louise Prinzessin in Preußen und Markgräfin zu Brandenburg, Äbtissin zu Herford*;[21,22] she was nineteen years old. From 1764 to 1802 she maintained her court in a manner befitting a royal household, administering the imperial convent and defending its rights.

Herford Abbey (*Frauenstift Herford*) was the oldest women's religious house in the Duchy of Saxony. Like other German convents, Herford was ruled by an Abbess with the title of Princess-Abbess (*Fürstäbtissin* or *Reichsäbtissin*). The abbey was dedicated in 832 and was elevated to the status of a Reichsabtei (Imperial abbey) under Emperor Louis the Pious (d. 840). In ecclesiastical matters, it answered directly to the Pope and was endowed with a third of the estates originally intended for Corvey Abbey. In 1533, during the Reformation, Herford Abbey became Lutheran, under the Electors of Brandenburg. From 1649 for over a century the abbesses were all Calvinist but that did not alter the Lutheran character of the principality.

When Friederike Charlotte became Abbess of Harford, she ranked among the independent princes of the Empire who lived in a princely state with a court of her own, and ruled their extensive convent estates. As Princess Abbess, Friederike Charlotte oversaw the principality of the Abby, which included about several thousand people, and the farms, vineyards, mills and factories. She could exercise supreme domestic authority over her monastery and all its dependencies, but as a female, she could not exercise any power of spiritual nature. However, the abbess was empowered to administer the temporal possessions of the convent; to issue commands to her nuns, binding them in conscience to maintain due discipline.

---

[20] Copy of *Lettres a une princesse d'allemagne* is available at www.e-rara.ch/doi/10.3931/e-rara-8748
[21] http://www.guide2womenleaders.com/womeninpower/Womeninpower1740.htm
[22] Herforder Fürstäbtissin Prinzessin Friederike Charlotte von Preußen (1764-1802)

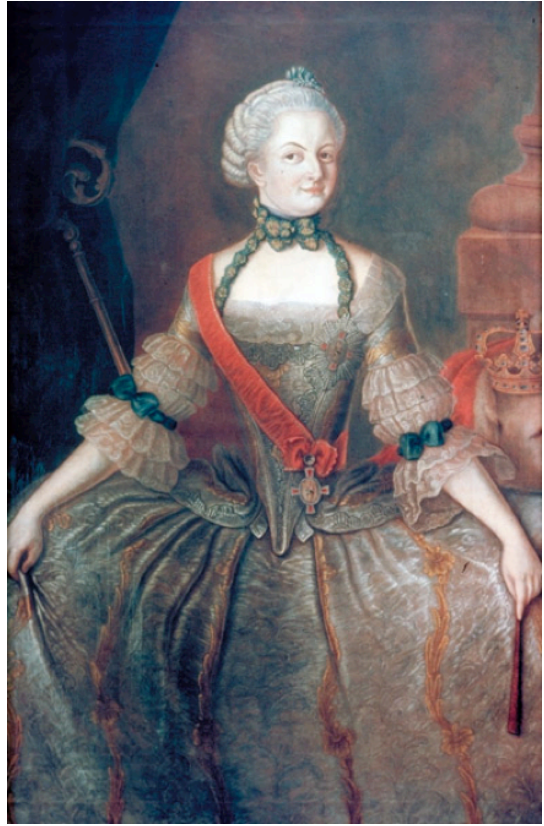

Friederica Charlotta Leopoldina Louise, Prinzessin in Preußen und
Markgräfin zu Brandenburg, Äbtissin zu Herford

We have no written testament describing the life and work of Friederike
Charlotte. We know that she was the last sovereign ruler of the Ecclesiastical Territory,
which was incorporated into Prussia in 1802 as part of the rearrangement of the German
Realm after the Napoleonic wars. On 15 August 1802, her Abbey was dissolved because
of secularization, and Frederick William III, King of Prussia, seized its assets. Princess
Friederike Charlotte remained with the collegiate ladies, who received a pension from the
kingdom. On 25 February 1803, Harford Abbey was annexed to Ravensberg within the
Kingdom of Prussia.

At that time, Europe was in turmoil. Kingdoms and nations were fighting for
increased power. In August 1805, Napoleon, Emperor of the French, turned his sights
from the English Channel to the Rhine in order to deal with the Austrian and Russian
threats. On 25 September, French troops invaded the Prussian realm. At the battle of
Austerlitz on 2 December 1805, Napoleon's army crushed the Third Coalition. At the
battle of Jena on 14 October 1806, French forces defeated the Prussian armies. With that,
Napoleon conquered almost all of Prussia, with the exception of Königsberg in East
Prussia (now it belongs to Russia). On 27 October 1806, Napoleon and his French Army
entered Berlin. By then, Frederick William III had fled with his court from Berlin to
Königsberg.

During the invasion, Friederike Charlotte, now sixty-one, fled to Altona (near Hamburg), where she died in 1808. She is buried in the church in Herford. The former abbey church remains in use today as Herford Minster (Herforder Münster).

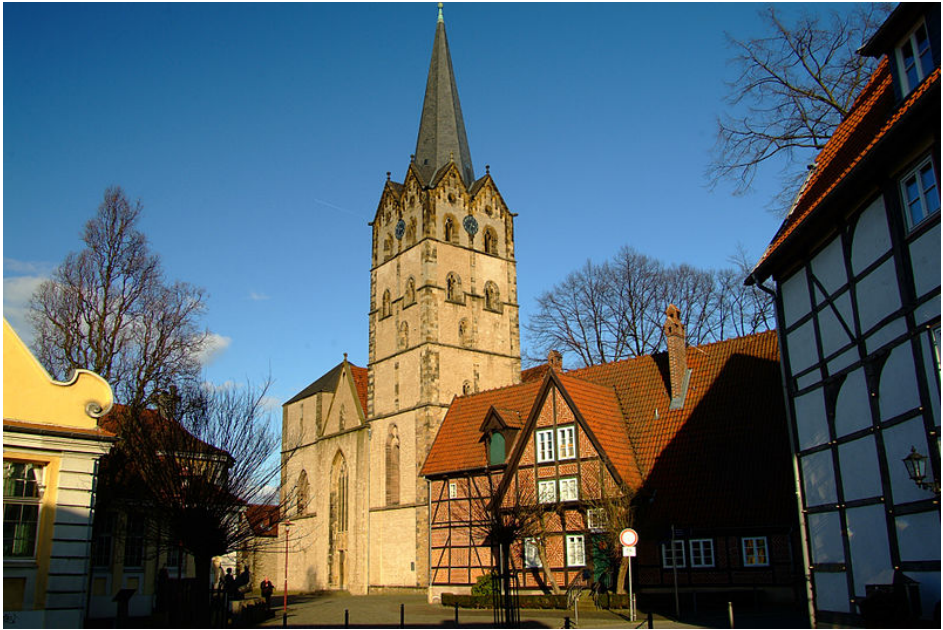

Former Herford Abbey church, now Herford Minster - During the middle ages it was a famous monastery for the education of aristocrat young women. The red building in the front is named the Kantorhaus. It was built around 1480, and it is the second oldest house in framework architecture in North Rhine-Westfalia.

## 5.    Afterthought

Why did Friederike Charlotte inherit the title Princess Abbess of Herford Abbey? Warfare was a significant factor in political competition, and marriage was one of the ways that alliances were made between the different kingdoms. Perhaps Friederike Charlotte did not have marriage suitors and thus becoming the Princess Abbess was a way to maintain the stronghold in that Prussian region. Her aristocratic familial connections and personal wealth must have been critical to the Abbey and Herford's economic viability.

In the eighteen-century, there were only two options for a girl of the nobility: marriage or a convent. The convent was a center of learning and also a refuge. Young ladies could live there if, for example, they had no dowry or their male relatives required all the family wealth. The convent could also serve as an asylum or a prison for women in cases of unfortunate marriages or injudicious love affairs. To the convent women could go if they had lost their beauty to the prevalent scourge of smallpox or had no marriage prospects.

A passing glance at Friederike Charlotte's portrait gives us some idea about her appearance. She seems to have been rather beautiful, and no historical account of a love interest exists. Her younger sister did get married at seventeen. On 25 July 1767, Louise Henriette Wilhelmine married her cousin Leopold III of Anhalt-Dessau in Charlottenburg (Berlin). By this marriage, Louise had the title of Princess, and later Duchess, of Anhalt-Dessau. We note that Euler never once mentioned the younger princess, suggesting that Louise Henriette did not receive his tutoring.

Friederike Charlotte was not the first princess to have been taught by a great scholar. Before her, in the seventeenth century, another woman of the nobility was linked through letters to a prominent mathematician and philosopher and she also became Abbess of Herford. Princess Elisabeth II of Bohemia (Pfalzgräfin bei Rhein or Countess Palatine of the Rhine) was an educated woman interested in the philosophy of René Descartes. Descartes wrote *La géométrie*, a book in which he included his application of algebra to geometry from which we now have Cartesian geometry. His work had a great influence on both mathematicians and philosophers. In 1642, the twenty-four year-old Princess Elisabeth read his *Meditations on First Philosophy*, and studied it thoroughly. A year later she wrote to Descartes expressing regret that they were unable to meet. She also posed some questions regarding his theory of the dualism of the body and soul. This began a long correspondence between Descartes and princess Elisabeth. The letters to each other have been preserved. At her request, Descartes became her teacher in philosophy and morals, and in 1644 he dedicated to her his *Principia Philosophiae*, a book where Descartes attempted to put the whole universe on a mathematical foundation, reducing the study to one of mechanics.

In 1649, Descartes accepted the invitation of Queen Christina of Sweden to go to Stockholm and be her tutor of philosophy. Descartes continued in correspondence with Elisabeth. However, their friendship would end rather quickly after his arrival to Sweden. Christina was a demanding monarch, expecting her lessons very early in the morning in a rather cold study. Descartes, who was accustomed to working in bed until noon, contracted pneumonia and died on 11 February 1650. Ten years after the death of Descartes, Princess Elisabeth entered Herford Abbey, and in 1667 she became Princess Abbess. During her reign, she distinguished herself by faithfulness in the performance of her duties, by her modesty and philanthropy, and especially by her kind hospitality. On February 8, 1680, Elizabeth died after a long and painful illness.

To date no manuscripts have been found about or written by the last Abbess of Herford Abbey. Thus, we don't know what Friederike Charlotte thought of Euler or what her achievements at Herford were. However, we can imagine her as a devout woman who, like her predecessors governed her Abbey with kindness and great wisdom. Of course, I wonder whether Friederike Charlotte kept Euler's original letters and if so, where are they now?

* * * * * * *